%% file: HyperSinger.tex
\documentclass[11pt]{amsart}

\usepackage{amssymb,amsmath,latexsym}
\usepackage{charter,tikz}
\usepackage{amssymb}
\usepackage{amscd}
\usepackage[all,cmtip]{xy}
\usepackage{rotating}
\usepackage{hyperref,mathrsfs}
\usepackage[osf]{mathpazo}

\newtheorem{theorem}{Theorem }[section]

\newtheorem{lemma}[theorem]{Lemma}

\newtheorem{remark}[theorem]{Remark}
\newtheorem{corollary}[theorem]{Corollary}
\newtheorem{proposition}[theorem]{Proposition}
\newtheorem{observation}[theorem]{Observation}

\newcommand{\bt}{\begin{theorem}}
\newcommand{\et}{\end{theorem}}
\newcommand{\bc}{\begin{corollary}}
\newcommand{\bl}{\begin{lemma}}
\newcommand{\ec}{\end{corollary}}
\newcommand{\el}{\end{lemma}}
\newcommand{\bo}{\begin{observation}}
\newcommand{\eo}{\end{observation}}
\newcommand{\bp}{\begin{proposition}}
\newcommand{\ep}{\end{proposition}}
\newcommand{\br}{\begin{remark}}
\newcommand{\er}{\end{remark}}
\newcommand{\brt}{\begin{result}}
\newcommand{\ert}{\end{result}}
\newcommand{\eop}{\hspace*{\fill}{\footnotesize $\blacksquare$}}

\def\I{\mathbf{I}}
\def\nI{\not\mathbf{I}}

\def\PG{\mathbf{PG}}
\def\PGL{\mathbf{PGL}}

\newcommand{\GL}{\mathbf{GL}}
\newcommand{\Aut}{\mathrm{Aut}}

\newcommand{\K}{\mathbb{K}}
\newcommand{\proj}{\mathrm{proj}}

\newcommand{\mF}{\mathscr{F}}
\newcommand{\hF}{\mathbf{F}}

\newcommand{\mP}{\mathscr{P}}
\newcommand{\mB}{\mathscr{B}}
\newcommand{\mL}{\mathscr{L}}

\newcommand{\mN}{\mathscr{N}}

\newcommand{\mS}{\mathscr{S}}
\newcommand{\cT}{\mathscr{T}}
\newcommand{\F}{\mathbb{F}}

\newcommand{\id}{\mathbf{1}}

\def\cO{\mathscr{O}}

\def\Z{\mathbb{Z}}

\def\F{\mathbb{F}}
\def\Fun{{\F_1}}

\def\Spec{\mathrm{Spec}}

\def\hF{\mathbf{F}}

\title[The Connes-Consani plane connection]{Hyperfield extensions, characteristic one and the Connes-Consani plane connection}


\author{Koen Thas}

\email{koen.thas@gmail.com}
\address{{Ghent University},
{Department of Mathematics},
{Krijgslaan 281, S25, B-9000 Ghent, Belgium}}
\urladdr{http://cage.ugent.be/$\sim$kthas}

\date{}

\begin{document}

\maketitle

\begin{abstract}
Inspired by a recent paper of Alain Connes and Catherina Consani which connects the geometric theory surrounding the elusive field with one element to sharply transitive group actions on finite and infinite projective spaces (``Singer actions''), we consider several fudamental problems and conjectures about Singer actions. Among other results, we show that virtually all infinite abelian groups and all (possibly infinitely generated) free groups act as Singer groups on certain projective planes, as a corollary of a general criterion. We investigate for which fields $\F$ the plane $\mathbf{P}^2(\F) = \PG(2,\F)$ (and more generally the space $\mathbf{P}^n(\F) = \PG(n,\F)$) admits a Singer group, and show, e.g., that for any prime $p$ and any positive integer $n > 1$, $\PG(n,\overline{\F_p})$ cannot admit Singer groups. One of the main results in characteristic $0$, also as a corollary of a criterion which applies to many other fields,  is that $\PG(m,\mathbb{R})$ with $m \ne 0$ a positive even integer, cannot admit Singer groups. 
\end{abstract}
\bigskip
{\footnotesize \tableofcontents}

\bigskip

\include{IntroAbsHyper}

\bigskip
\include{F1Hyper}

\bigskip
\vspace*{1cm}\section{What is known}

Karzel proves the following (answering a more general version of a question of Hall \cite{Hall}):

\bt[H. Karzel \cite{Karzel}]
\label{Karz}
Let $S$ be a finitely generated commutative Singer group of $\PG(m,\F)$, with $m \in \mathbb{N}^\times$ and different from $1$,
and $\F$ a field. Then $\F$ is finite, and $S$ is cyclic.
\et

So Desarguesian projective spaces (different from projective lines) can only allow  commutative Singer groups which are 
infinitely generated. Later, we will construct commutative Singer groups with this property for spaces $\PG(m,\F)$ for many values of $(m,\F)$.

On the other hand, we will show in the next section that virtually any infinite  commutative group (those that do not have involutions) can 
act as a Singer group on an appropriate projective plane | so also the finitely generated examples | but the planes are not 
Desarguesian by Karzel's result. Our result is a corollary of a theorem which generalizes the next result of Hughes:

\bt[D. R. Hughes \cite{Hughes}]
\label{Hugh}
Let $H$ be a countably infinite group, and
assume the following properties for $H$:
\begin{equation}
\left\{\begin{array}{ccccc}
(h_1) &h^2 &\ne & \id &\forall h \in H^\times \\
(h_2) &\vert h^H \vert &= &\infty &\forall h \in H \setminus Z(H) \\
(h_3) &\#\{ x \vert x^2 = h'\}&<&\infty &\forall h' \in H \\
\end{array}\right.
\end{equation}
Here, $Z(H)$ is the center of $H$. Then $H$ acts as a Singer group on some projective plane.
\et

Hughes applied Theorem \ref{Hugh} to show that free groups with a finite number $n$ of generators ($n \geq 2$) act as Singer groups of 
some projective plane. Later, we will obtain this result for any free group.\\

For vector spaces, one could also consider the related problem of studying sharply transitive automorphism groups of the nonzero vectors. 
We recall the following nice result.

\bt[\cite{Cheretal}]
Let $\mathbb{K}$ be an algebraically closed field and $G$ a subgroup of $\mathbf{GL}_n(\mathbb{K})$ which acts sharply transitively on the set of nonzero vectors in $\mathbb{K}^n$, where $n \in \mathbb{N}^\times$. Then either $n = 1$ and $G = \mathbb{K}^\times$, or $n = 2$, and $G$ can be precisely described.
\et 

The vector space problem relates to our problem as follows.
Let $V$ be an $\F$-vector space over the field $\F$, and let $\PG(V)$ be the corresponding projective space. If $K$ is a Singer group of $\PG(V)$, then if $\overline{K}$
is the corresponding automorphism group of $\PG(V)$ (that is, $\overline{K}$ contains the full scalar group $S$ and $\overline{K}/S = K$), the latter acts 
sharply transitively on the nonzero vectors. And if $K$ is linear, $\overline{K}$ is linear as well. (Of course, if $\F$ is algebraically closed, Singer groups cannot exist
due to the fact that all polynomials over $\F$ have roots in $\F$, cf. later sections.) Conversely, if $H$ acts sharply transitively on the nonzero vectors of 
the vector space $V$, and $H$ is abelian, then $H$ contains all scalar automorphisms $S$\footnote{If an element of $H$ fixes some vector line, it must fix all vector lines as $H$ is abelian and transitive on vector lines, so due to the transitivity on nonzero vectors, $H$ contains all scalar automorphisms.}, and $H/S$ induces an abelian Singer group of the associated projective space.

As we will construct Singer groups for ``most'' projective spaces, we will get the sharply transitive groups of the vector space for free.

\bigskip
\vspace*{1cm}\section{Construction of Singer groups}

In this section, we slightly generalize the result of Hughes on planar difference sets in not necessarily abelian groups,
by removing the assumption on countability. \\

\subsection{Partial difference sets}

If $G$ is a group (written multiplicatively)  and $C \subseteq G$, denote by $D(C)$ the set of ``differences'' $\{ c{c'}^{-1} \vert c, c' \in C\}$. Now assume that $K$ is a group, 
and $S$ a subset such that for any $k \in K^\times$, there is precisely one couple $(a,b) \in S \times S$ such that $k = ab^{-1}$ | in other words, the map
\begin{equation}
\phi: S \times S \setminus \mathrm{diagonal} \longrightarrow K^\times: (a,b) \longrightarrow ab^{-1}
\end{equation}
is a bijection (and as a consequence, $D(C) = K$). We call $S$ a {\em difference set} in $G$. If the map $\phi$ merely is injective, we call $S$ a {\em partial difference set}.
Then defining a ``point set'' $\mP = G$ and ``line set'' $\mB = G$, where a point $x$ is incident with a line $y$ (and we write $x \I y$) if and only if $xy^{-1} \in S$, we 
obtain a projective plane $\Gamma(G,S)$ with the special property that  $G$ acts (by right translation) as a sharply point transitive automorphism group (= Singer group).\footnote{For $g \in G$, we have $a \I b$ if and only if $ab^{-1} \in S$ if and only if $(ag)(g^{-1}b^{-1}) \in S$ if and only if $ag \I bg$.} And
conversely, a projective plane admitting a Singer group $G$ can always be constructed in this way from a difference set $S \subset G$.

We will use the following easy lemma without reference.
\bl
If ${(S_\omega)}_{\omega \in \Omega}$ is a chain of partial difference sets in a group $K$ (so $\Omega$ is well-ordered and from 
$\nu < \mu$ follows that $S_{\nu} \subseteq S_{\mu}$), then $\cup_{\omega \in \Omega}S_{\omega}$ also is a partial difference set.\\ 
\el

\medskip
\subsection{Ordinals}

Each {\em ordinal} is the well-ordered set of all smaller ordinals. The smallest ordinal is $0 = \emptyset$,  the next-smallest ordinal is
$1 = \{0\} = \{\emptyset\}$, followed by
$2 = \{0,1\} = \{\emptyset,\{\emptyset\}\}$, etc. After all finite ordinals have been constructed, we continue with
$\omega = \{0,1, 2,\ldots \}$, $\omega + 1 = \{0,1,2,\ldots\} \cup \{\omega\}$, and so on,
and eventually $2\omega = \{0,1,2,\ldots\} \cup \{\omega,\omega + 1,\omega + 2,\ldots \}$.
This is followed by $2\omega + 1,2\omega + 2, \ldots, \omega^2$. And so on. All the ordinals we have mentioned so far are countably infinite. After all countable ordinals have been defined, we
meet the first uncountable ordinal, denoted $\omega_1$; later we reach $\omega_2$, etc. Let $\gamma$ be an ordinal. 
Then the {\em successor} of $\gamma$ is 
\begin{equation}
\gamma + 1 =  \gamma \cup \{\gamma\}, 
\end{equation}
this being the smallest ordinal exceeding $\gamma$. Every ordinal is either a successor ordinal or a {\em limit ordinal}, but never both. 
A {\em limit ordinal} is an ordinal $\alpha$ such that 
\begin{equation}
\alpha = \bigcup_{\beta<\alpha}\beta.
\end{equation}
Now an arbitrary set $S$ may be indexed as $S = \{ s_{a} \vert a \in A\}$,
where $A$ is an ordinal. Moreover we may assume $A$ to be {\em minimal} among all
ordinals of cardinality $\vert A \vert$ | otherwise we may simply re-index suitably. 

\medskip
\subsection{Construction}

Let $H$ be an infinite group, and let $\vert H\vert = A$ be the smallest ordinal of cardinality $\vert H\vert$;
write $H = \{ h_{\alpha} \vert \alpha \in A \}$ ($A$ is well-ordened).

Define for each $\gamma \in A$ a set $S_{\gamma}$ such that
\begin{itemize}
\item[(i)]
$\vert S_{\gamma}\vert \leq \vert\gamma\vert < \vert A\vert$;
\item[(ii)]
$h_{\gamma} \in D(S_{\gamma})$ for $\gamma \in A$;
\item[(iii)]
$S_{\gamma}$ is a partial difference set of $H$;
\item[(iv)]
$S_{\gamma} \subseteq S_{\beta}$ for $\gamma < \beta$ and $\beta \in A$.
\end{itemize}

If $\alpha$ is a limit ordinal, define $S_{\alpha} = \cup_{\beta < \alpha}S_{\beta}$. (Note that $\vert S_{\alpha}\vert \leq \vert \alpha \vert \cdot \vert \alpha \vert$.) 
Now let $\alpha$ be a successor ordinal $\alpha = \beta + 1$; if $h_{\alpha} \in D(S_{\beta})$, put $S_{\alpha} = S_{\beta}$.
Otherwise, we construct $S_{\alpha}$ by adding two new elements to $S_{\beta}$ such that $h_{\alpha} \in D(S_{\alpha})$.

We seek properties for $H$ such that this particular step (and then the whole construction) can be carried out.
Let $S_{\beta} = \{ s_i \vert i \in I \}$ with $\vert I \vert < \vert A \vert$. (We suppose without loss of generality that $\vert I \vert \not\in \mathbb{N}$.) 
Suppose $d = h_{\alpha} \not\in D(S_{\beta})$.

Assume the following properties for $H$:
\begin{equation}
\left\{\begin{array}{ccccc}
(d_1) &h^2 &\ne & \id &\forall h \in H^\times \\
(d_2) &\vert h^H \vert &= &\vert H \vert &\forall h \in H \setminus Z(H) \\
(d_3) &\#\{ x \vert x^2 = h'\}&<&\vert H\vert &\forall h' \in H \\
\end{array}\right.
\end{equation}

Note that a planar  Singer group never can have involutions (if $\sigma$ would be such an involution and $L$ is a line of the plane,
$L \cap L^{\sigma}$ would be a fixed point), so the first property is necessary.\\

Take any $x \in H \setminus S_{\beta}$, and consider the following sets

\begin{equation}
\left\{\begin{array}{c}
D(S_{\beta} \cup \{x\})\\
D(S_{\beta} \cup \{d^{-1}x\})\\
D(S_{\beta})\\
\end{array}\right.
\end{equation}

(Note that $\vert D(S_{\beta})\vert = \vert S_{\beta}\vert$.) We need to check that there exist $x$ such that 
$S_{\beta} \cup \{ x,d^{-1}x\}$ is a partial difference set, observing that any couple $(a,b)$ for which $ab^{-1} = d$ indeed has the form 
$(x,d^{-1}x)$. 

First note that, with $s_i, s_j \in S_{\beta}$, if $xs_i^{-1} = d^{-1}xs_j^{-1}$, then $d^x = s_j^{-1}s_i$, and if $xs_i^{-1} = s_jx^{-1}d$, then 
$s_j^{-1}ds_i = (s_j^{-1}x)^2$. Besides that, the number of elements $x$ for which an element of $D(S_{\beta} \cup \{x\}) \cup
D(S_{\beta} \cup \{d^{-1}x\})$ equals an element of $D(S_{\beta})$ is at most $\vert S_{\beta}\vert$. So ($d_2$) and ($d_3$) readily imply 
that we can find $x$ such that $S_{\beta} \cup \{ x,d^{-1}x\}$ is a partial difference set. 

We have proved the following theorem:

\bt[Construction]
If an infinite group $H$ satisfies the following properties, then $H$ acts as a Singer group on some projective plane.
\begin{equation}
\left\{\begin{array}{ccccc}
(d_1) &h^2 &\ne & \id &\forall h \in H^\times \\
(d_2) &\vert h^H \vert &= &\vert H \vert &\forall h \in H \setminus Z(H) \\
(d_3) &\#\{ x \vert x^2 = h'\}&<&\vert H\vert &\forall h' \in H \\
\end{array}\right.
\end{equation}
\eop
\et

\medskip
\subsection{An example: general free groups}

In \cite{Hughes}, Hughes showed that free groups on a finite number of generators satisfy the properties (h$_1$)-(h$_2$)-(h$_3$) of Theorem \ref{Hugh}, 
so that they act on certain projective planes as Singer groups. Now let $\hF(\Omega)$ be a free group with generator set $\Omega$, where $\Omega$
is any infinite alphabet. For any reduced element $f \in \hF(\Omega)$, let $\pi(f)$ be the subset of $\Omega \cup \Omega^{-1}$ of letters used in $f$. 

First of all, note that $\hF(\Omega)$ cannot have involutions since if $x$ would be an involution, it also would be an involution in $\hF(\pi(x))$, which 
contradicts Hughes's result.

Next, let $h$ be any nontrivial element in $\hF(\Omega)$, and consider the equation
\begin{equation}
x^2 = h.
\end{equation}
Hughes shows in \cite{Hughes} that this equation has a unique solution in a free group $\hF(S)$ where $\pi(h) \subseteq S$ and $S$ is finite, so
it follows easily that it also has a unique solution in $\hF(\Omega)$.

Next we want to consider orbits $g^{\hF(\Omega)}$. Define $\Omega' := \Omega \setminus \pi(g)$, and define the set
\begin{equation}
\xi(g) := \{g^{\omega}:= \omega^{-1}g\omega \vert \omega \in \Omega'\} \subset g^{\hF(\Omega)}.
\end{equation}
It follows that
\begin{equation}
\vert \xi(g)\vert = \vert \Omega' \vert = \vert \Omega\vert = \vert \hF(\Omega)\vert.
\end{equation}

So indeed (d$_1$)-(d$_2$)-(d$_3$) are satisfied, and whence $\hF(\Omega)$ acts as a Singer group on some plane.

\bigskip
\vspace*{1cm}\section{Construction of difference sets | Abelian case}

By ($d_2$), one would expect that the previous section would not apply to the abelian case, but this is, in fact, not entirely true.

For suppose $H$ is abelian now, without involutions (cf. ($d_1$)).
If, as above, we want to add $x$ and $d^{-1}x$ to $S_{\beta}$ to obtain a partial difference set $S_{\beta} \cup \{ x,d^{-1}x\}$ for which $d$
is a difference, we need to find an $x$ for which $d^x \ne s_j^{-1}s_i$; but from $d^x = {s_j}^{-1}s_i$ we would have $d^x = d = s_j^{-1}s_i = s_is_j^{-1}$ since $H$ is abelian, contradiction
since $d \not\in D(S_{\beta})$. So ($d_2$) is not needed here. Secondly, suppose 
\begin{equation}
\#\{ x \not\in S_{\beta} \vert  s_j^{-1}ds_i = (s_j^{-1}x)^2\ \mbox{for some}\ s_i,s_j \in S_{\beta}\} > \vert S_{\beta} \vert.
\end{equation}
Then obviously we can find an $s_\ell \in S_{\beta}$ and different $z, z' \not\in S_{\beta}$ such that $(s_{\ell}^{-1}z)^2 = (s_{\ell}^{-1}{z'})^2$,
implying that $z{z'}^{-1}$ is an involution, contradiction. So for abelian groups, ($d_3$) need not be assumed since it follows (in the context that we need it) from
($d_1$).

So for the abelian case, we obtain the most general constructive result as possible:

\bt[Abelian Singer groups | Characterization]
An infinite abelian group acts as a Singer group on some projective plane if and only if it contains no involutions.
\eop
\et

\bigskip
\vspace*{1cm}\section{Singer groups for classical spaces}
\label{class}

Let $\F$ be any field,  and suppose $\F'/\F$ is a proper  field extension. Then $\F'$ can be naturally seen
as an $\F$-vector space $V(\F')$ as before. Now ${\F'}^{\times}$ acts by multiplication on $V(\F')$, and clearly this induces a subgroup of $\GL(V(\F'))$ which
acts sharply transitively on the nonzero vectors. The subgroup $\F^{\times} \leq {\F'}^{\times}$ acts as scalars, and ${\F'}^{\times}/\F^{\times}$ induces 
a sharply transitive group on the points of the projective space $\PG(V(\F'))$. 
If the degree of $\F'/\F$ is a natural nonzero number $m$ (which is at least $2$), then $\PG(V(\F')) \cong \PG(m - 1,\F)$. We have

\bt
\label{6.1}
If $\omega = [\F' : \F]$ is the not necessarily finite degree of the field extension $\F'/\F$, then $\PG(\omega - 1,\F)$ allows a linear Singer  group.
In particular, if $\omega = 3$, this applies to the Desarguesian plane $\PG(2,\F)$.
\eop \\
\et

Theorem \ref{6.1} motivates us to state the next simple idea as a principle.\\
 
\quad\textsc{Singer|Algebraic closure principle}.\quad
{\em The farther away a field $\F$ is from its algebraic closure, the more Desarguesian projective spaces over $\F$ allow a (linear) Singer group in this 
construction. And the more isomorphism classes of such groups arise.}\\

For $n \in \mathbb{N} \setminus \{0,1\}$, call a field $\F$ {\em $n$-ally closed} if every polynomial of degree $n$ in $\F[x]$ as a root.
When $n = 2$, we also speak of ``quadratically closed,'' when $n = 3$ ``cubically closed.''

\bc
If a field $\F$ is not $n$-ally closed, $n \in \mathbb{N} \setminus \{0,1\}$, then $\PG(n - 1,\F)$ admits a (linear) Singer group.
In particular, if a field $\F$ is not cubically closed, then $\PG(2,\F)$ admits a (linear) Singer group.
\eop \\
\ec

The next corollary is a nice example of the aforementioned pinciple.

\bc
If a field $\F$ is not algebraically closed, then $\PG(\ell,\F)$ admits a (linear) Singer group for some $\ell \in \mathbb{N} \setminus \{0,1\}$.
\ec
{\em Proof}.\quad
As $\F$ is not algebraically closed, there exists an element in $\F[x]$, say of degree $\ell > 1$, without roots. Now apply Theorem \ref{6.1}.
\eop \\

(In the next section, one can find more formal information about real-closed fields.)

\bigskip
\vspace*{1cm}\section{$N$-Ally closed fields with few automorphisms}

For some fields $\F$, it is rather easy to exclude the existence of Singer groups for $\PG(n - 1,\F)$, $n \in \mathbb{N} \setminus \{0,1,2\}$, $n$ arbitrary odd. 

\bt
\label{struct}
Suppose $\overline{\F}$ is such that $[\overline{\F} : \F]$ is finite of degree $m \ne 1$. Suppose furthermore that 
\begin{equation}
\vert \Aut(\F) \vert < \vert \F\vert 
\end{equation}
(where $\vert \cdot \vert$ denotes cardinality and $\Aut(\cdot)$ the automorphism group).  Then $\PG(n - 1,\F)$ does not admit a Singer group, where  $n \in \mathbb{N} \setminus \{0,1,2\}$, $n$ odd.
\et
{\em Proof}.\quad
Suppose by way of contradiction that $S$ is a Singer group of $\PG(n - 1,\F)$. Let $\mP$ be the point set of $\PG(n - 1,\F)$. As $\vert S \vert = \vert \mP\vert = 
\vert \F\vert > \vert \Aut(\F)\vert$, $S$ must have a nontrivial intersection with $\GL_{n}(\F)$. For, the map
\begin{equation}
\pi: S \longrightarrow \Aut(\F): (A,\sigma) \longrightarrow \sigma,
\end{equation}
cannot be injective, so that there are different elements of the form $(B,\kappa)$, $(C,\kappa)$ in $S$; appropriately composing yields nontrivial linear
elements in $S$.
So suppose $\gamma \in S^\times$ is a linear element, 
corresponding to the matrix $A \in \GL_n(\F)$. 
Fields $\F$ for which $[\overline{\F} : \F]$ is finite are classified by a result of Artin-Schreier \cite{ArtSchr}: either $\F$ is algebraically closed (so the degree is $1$), or 
$[\overline{\F} : \F] = 2$ and $\F$ is real-closed (and conversely).
As $(n,[\overline{\F} : \F]) = 1$, the characteristic polynomial of $A$ must have at least one 
root $\rho$ in $\F$. If $v$ is a $\rho$-eigenvector, the corresponding point of $\PG(n - 1,\F)$ is fixed by $\gamma$, contradiction. \eop \\

(In the above, it makes no sense to allow the extension $[\overline{\F} : \F] = 1$, since $\vert \Aut(\overline{\F})\vert = \#2^{\overline{\F}}$.)

\bc
If $\F$ is real-closed and the positive odd integer $n$ is at least $3$, $\PG(n - 1,\F)$ cannot admit Singer groups if $\vert \Aut(\F) \vert < \vert \F\vert$. In particular, 
$\PG(2,\mathbb{R})$ has no Singer groups.
\ec
{\em Proof}.\quad
We have that $\vert \Aut(\mathbb{R}) \vert = 1$. \eop \\

Now let $\F$ be a real-closed field. Then for any $k \in \F$, either $k$ or $-k$ is in $\F^2$ (the set of squares). There is a unique 
total order $\leq$ on $\F$ defined by:
\begin{equation}
0 \leq y \ \mbox{if and only if}\ y \in \F^2.
\end{equation}
(The order is unique as squares must be positive with respect to a total order.) Let $\alpha \in \Aut(\F)$; as $\alpha(k^2) = \alpha(k)^2$ for any
$k \in \F$, $\alpha$ preserves the order (as $a < b$ if and only if there is a nonzero square $c^2$ such that $a + c^2 = b$).

The next result detects certain real-closed fields with trivial automorphism groups.

\bt
Let $\F$ be a real-closed field which is a subfield of $\mathbb{R}$. Then $\Aut(\F)$ is trivial.
\et
{\em Proof}.\quad
Let $\beta \in \Aut(\F)$; then $\mathbb{Q} \leq \F$ is fixed elementwise, and $\beta$ preserves the unique total order $\leq$ on $\F$ (which is 
the one inherited by $\mathbb{R}$). For $\kappa \in \F$, define $\mathbb{Q}^+(\kappa) := \{ q \in \mathbb{Q} \vert q \geq \kappa \}$ and 
$\mathbb{Q}^-(\kappa) := \{ q \in \mathbb{Q} \vert q \leq \kappa \}$. Note that both $\mathbb{Q}^+(\kappa)$ and $\mathbb{Q}^-(\kappa)$ are uniquely defined
by $\kappa$; if $\kappa \ne \kappa'$ are elements of $\F$ and $\kappa < \kappa'$, then there is a rational number $q$ such that $\kappa < q < \kappa'$.
Whence $q \in \mathbb{Q}^+(\kappa) \cap \mathbb{Q}^-(\kappa')$. It follows that all $\kappa \in \F$ must be fixed by $\beta$.
\eop \\

\bc
\label{exmp}
 For the positive odd integer $n$ which is at least $3$, $\PG(n - 1,\F)$ cannot admit Singer groups if $\F$ is either the field of real algebraic numbers, the field of computable numbers or the field of (real) definable numbers. In particular, 
$\PG(2,\F)$ has no Singer groups in these cases.
\ec
{\em Proof}.\quad
Each of these fields is real-closed and a subfield of the reals. By the previous theorem, $\Aut(\F)$ is always trivial. The statement then follows from Theorem \ref{struct}.
\eop \\

Isolating the dimension, Theorem \ref{struct} generalizes immediately as follows.

\bt
\label{structgen}
Let $\F$ be $\ell$-ally closed with $\ell \in \mathbb{N} \setminus \{ 0,1\}$. Suppose furthermore that 
\begin{equation}
\vert \Aut(\F) \vert < \vert \F\vert. 
\end{equation}
Then $\PG(\ell - 1,\F)$ does not admit a Singer group.
\eop \\
\et

\bigskip
\vspace*{1cm}\section{Possible strategy for classification?}
\label{strat}

One is tempted to study the following property (which we formulate for planes, but which is easily generalized to other spaces):\\

\quad(E)\quad {\em Let $\F$ be a field and $\mathbb{K}\vert\F$ a field extension. 
If $S$ is a Singer group of $\PG(2,\F)$, then $\PG(2,\mathbb{K})$ also allows a Singer group  $S'$ such that $S'_{\vert \F} = S$.}\\

In a category $\mathbf{E}$ of fields for which (E) is true (completed by the appropriate fields), we have:\\

\quad(AC)\quad {\em If $\F$ is an object in $\mathbf{E}$, then $\PG(2,\cup_{\F' \in \mathbf{E}, \F \leq \F'} \F')$ also allows a Singer group.}\\

(Consider an arbitrary filtration
\begin{equation}
\F = \F^{(0)} \leq \F^{(1)} \leq \cdots 
\end{equation}
such that $\cup_{i \in \mathbb{N}}\F^{(i)} = \cup_{\F \leq \F' \in \mathbf{E}}\F'$, and take a direct limit of the induced directed system of Singer groups.)

Let us for instance define a category $\mathbf{E}$ as having as objects a fixed finite field $\F_p$, $p$ a prime, and all finite extensions $\F_{p^i}$ with
$(i,3) = 1$. (Morphisms are natural.) Then by \S \ref{class}, we can construct a canonical Singer group $S(\PG(2,q))$ for each object $\F_q$ in $\mathbf{E}$.
As we will later see, the property $(3,i) = 1$ translates in the fact that if $m$ divides $n$, $m, n \in \mathbb{N} \setminus 3\mathbb{N}$, then 
$S(\PG(2,p^m)) \leq S(\PG(2,p^n))$. Taking the direct limit of the naturally defined directed system of groups, we obtain a Singer group of 
$\PG(2,\cup_{\F \in \mathbf{E}}\F)$. Of course, this specific Singer group can also be obtained directly by using \S \ref{class} (since this 
limit has field extensions of degree $3$). \\

If (E) would be true for a sufficiently large category of field extensions of a fixed field $\K$, and $\overline{\K}$ is an algebraically closed field for which $\PG(2,\overline{\K})$ does not have a Singer group, then 
$\PG(2,\K)$ also does not have a Singer group. Unfortunately, even for the category for finite fields (with completions)
(E) is not satisfied, although almost.  Still, all categories $\mathbf{E}$ of field extensions of some fixed field $\F$ that enjoy (E) also enjoy (AC), so that 
Singer groups for the ``$\mathbf{E}$-closures'' also exist.

\bigskip
\vspace*{1cm}\section{Algebraically closed fields}

In the case of algebraically closed fields, we can say the following.

\bt
\label{torsion}
Let $S$ be a Singer group of $\PG(m,\F)$, $m \in \mathbb{N}^{\times}$, $\F$ algebraically closed. Then 
$S$ is torsion-free if $\mathrm{char}(\F) \ne 0$. If $\mathrm{char}(\F) = 0$, then $S$ is torsion-free if $m = 2$.  
\et
{\em Proof}.\quad
Let us represent an element $\mu$ of $\mathbf{P\Gamma L}_{m + 1}(\F)$ by $(A,\sigma)$, where $A$ is the corresponding element in $\GL_{m + 1}(\F)$, and 
$\sigma \in \Aut(\F)$, and where the automorphism is defined as $x \longrightarrow Ax^{\sigma}$ (in homogeneous coordinates with respect to some 
coordinate system). Suppose there is some $\ell \in \mathbb{N}^\times$ such that $(A,\sigma)^{\ell} = (\lambda\mathbf{I}_{m + 1},\id)$, $\lambda \ne 0$ (and assume $\ell$ is the order of $\mu$ as an element of $S$). 
Denoting 
by $\F_{\sigma}$ the fixed field of $\sigma$, we have that $\ell \geq \vert \langle \sigma \rangle \vert = [\F : \F_{\sigma}] < \infty$ ($\F/\F_{\sigma}$ is Galois by a theorem of Artin), so that by the Artin-Schreier Theorem \cite{ArtSchr},  $\vert \langle \sigma \rangle \vert \leq 2$. 

When the characteristic is not zero, $\Aut(\F)$ is torsion-free by e.g. a result of Baer \cite{Baer}, so that necessarily $\sigma = \id$.
As before, it follows that $\mu$ has fixed points.

Let the characteristic be zero. Clearly we must assume that $\sigma$ is a nontrivial involution. If $\ell > 2$, then $(A,\sigma)^2 = (AA^\sigma,\id)$ induces
a nontrivial linear element of $S$, and then we know it has fixed points. So $\mu$ is an involution.

Finally, Singer groups cannot have involutions in the planar case (in any characteristic).
\eop \\

\br
{\rm Note that by the proof of the previous theorem, the projection of $S$ onto $\Aut(\F)$ is also torsion-free, as is its projection onto $\PGL_{m + 1}(\F)$.}
\er

\br
{\rm Note that if an involution $\sigma$ of some projective space in characteristic $0$ does not fix any point, it must fix a parallel class of lines elementwise.}
\er

\medskip
\subsection{Singer groups of $\PG(2,\overline{\mathbb{F}_p})$ do not exist}
\label{finac}

Put $\mathbb{N}^\times =: I$, and make the latter into
a directed set, by writing that $n \preceq m$ if $n \vert m$ and $(m/n,3) = 1$. We will first explain the motivation for this definition.  
Let $p$ be a prime, and consider $\F_i := \F_{p^i} \leq \F_j  := \F_{p_j}$ with $i \preceq j \ne i$. Let $\F_j' = \F_j[X]/(f(X))$ be an extension of degree $3$ of $\F_j$
($f(X)$ having degree $3$), and define $\F_i' := \F_i[X]/(f(X))$ | this extension is also of degree $3$, and is a subfield of $\F_j'$. Then $\F_j \cap \F_i' = \F_i$ as 
$(j/i,3) = 1$. We have that 

\begin{equation}
\label{eqcalc}
{\F_i'}^{\times}/\F_i^{\times} = {\F_i'}^{\times}/({\F_i'}^{\times} \cap \F_j^{\times}) \cong  {\F_i'}^{\times}{\F_j^{\times}}/\F_j^{\times} \leq {\F_j'}^{\times}/\F_j^{\times}.
\end{equation}

In other words, we have an inclusion of cyclic groups $C_{p^{2i} + p^i + 1} \leq C_{p^{2j} + p^j + 1}$.\footnote{Notice that this formula gives an easy, calculation-free, proof of the following property:
\begin{lemma}
If $j \equiv 0\mod{i}$, $i$ and $j$ being positive nonzero integers, and $(j/i,3) = 1$, then $p^{2i} + p^i + 1$ divides $p^{2j} + p^j + 1$ for any prime $p$.
\eop
\end{lemma}}

This is exactly what we need for having a good definition for the directed system above | for the infinite case, we will use the form of $(\ref{eqcalc})$.
Unfortunately, our system is not directed anymore due to the divisibility constraint: if $i, j \in I$ and $3^n \vert i$ but not $3^{n + 1}$, and $3^m \vert j$ but not 
$3^{m + 1}$ and $n \ne m$, then there is no $k \in I$ such that $i \preceq k$ and $j \preceq k$. (Similar obstructions arise when going to higher dimensions.)

The next theorem explains that we can not adapt the construction.

\bt[Nonexistence for the fields $\overline{\F_p}$]
For any prime $p$ and any positive integer $m \geq 1$, the space $\PG(m,\overline{\F_p})$ does not admit Singer groups. In particular, this result 
applies to the planar case.
\et
{\em Proof}.\quad
Let $p, m$ and $\PG(m,\overline{\F_p})$ as in the statement, and suppose that $S$ is a Singer group. We represent an element $\gamma$ of $S$, as before, by a couple $(A,\sigma)$, where $A \in \GL_{m + 1}(\overline{\F_p})$ and $\sigma \in \Aut(\overline{\F_p})$. If $\sigma = \id$, we know that $\gamma$ has fixed points, so $\sigma \ne \id$ and by Theorem \ref{torsion} we have that $\langle \gamma \rangle \cong \mathbb{Z}, +$. Write $\overline{\F_p}$ as $\cup_{i = 1}^{\infty}\F_{p^i}$. As $A$ has a finite number of entries, and 
as each element of $\overline{\F_p}$ is contained in some finite subfield, there is a finite subfield $\F_q$ such that  $A \in \GL_{m + 1}(q)$. For each $i \in \mathbb{N}^\times$, $\overline{\F_p}$ contains a unique subfield of size $\F_{p^i}$, so $\sigma$ stabilizes all these subfields | in particular, it stabilizes $\F_q$. So $\gamma$ fixes $\PG(m,q) \subseteq
\PG(m,\overline{\F_p})$. As $\langle \gamma \rangle$ is not finite, some power of $\gamma$ fixes points (and even all points) of $\PG(m,q)$, contradiction. (Another way of finishing the proof, knowing that $A \in \GL_{m + 1}(\F_q)$, is to use the remark after Theorem \ref{torsion}.)
\eop \\

\subsection{Singer actions for $\PG(r,\overline{\mathbb{Q}})$}

Recalling the result of Karzel stated in the beginning of this paper, we deduce the following now.

\bt
\label{strfree}
Let $\F$ be a field which is not finite.
Let $S$ be a Singer group of $\PG(r,\F)$, $r \in \mathbb{N} \setminus \{ 0,1\}$,
and suppose $X \subseteq S$ is an independent set of generators  of $S$ of minimal size $\omega$.
Then either $\omega = \vert \F\vert$, or $\F$ is countably infinite. If in the latter case $w \ne \vert \F\vert$ (so if $\omega$ is finite), $\F$ is not isomorphic to $\overline{\mathbb{Q}}$.
\et

{\em Proof}.\quad
As the free group $\mathbf{F}_\omega$ on $\omega$ letters has size $\mathrm{max}(w,\vert \mathbb{Q} \vert)$, and as each group generated by $\omega$
elements is a homomorphic image of $\mathbf{F}_{\omega}$, it is clear that we only have to consider the case where $\F$ is countable. 
(If $\F$ is not countable, $S$ is not countable. As $S \cong \mathbf{F}_{\omega}/\mbox{kernel}$, it follows that $\omega$ is not countable. So 
$\omega = \vert \F\vert$.)
So we suppose that $\omega \in \mathbb{N}$, and that $\F$ is isomorphic to $\overline{\mathbb{Q}}$ by way of contradiction. Suppose $(A_1,\sigma_1),\ldots,(A_\omega,\sigma_\omega)$ represent the elements of $X$ ($A_i \in \GL_3(\F)$, $\sigma_i \in \Aut(\F)$).
If $\Omega$ is the set of entries of $A_1,\ldots,A_{\omega}$, $\Omega$ is a finite set of algebraic numbers. So there is an algebraic number $\alpha$ such that 
the field extension $\mathbb{Q}(\Omega)$ coincides with $\mathbb{Q}(\alpha)$.  Suppose $[\mathbb{Q}(\alpha) : \mathbb{Q}] = m \in \mathbb{N}^\times$. 
Then $\vert \mathbb{Q}(\alpha)^{\Aut(\overline{\mathbb{Q}})} \vert = m$ (as the precise number of monomorphisms $\gamma: \mathbb{Q}(\alpha) \longrightarrow \mathbb{C}$
is given by $m$, essentially by sending $\alpha$ to its conjugates over $\mathbb{Q}$). It follows that $\mathbb{Q}(\alpha)^\Sigma$, where $\Sigma = \langle \sigma_i \vert i = 1,\ldots,\omega \rangle$,
generates a proper subfield $\F'$ of $\overline{\mathbb{Q}}$ (which is finitely generated over $\mathbb{Q}$). Moreover, the space $\PG(r,\F')$ which is a proper
subspace of $\PG(r,\overline{\mathbb{Q}})$ by field reduction, is stabilized by each element of $X$, so also by $S$. But this contradicts the fact that $S$ is transitive on the points of $\PG(r,\overline{\mathbb{Q}})$.
\eop

\bc
Let $S$ be a Singer group of $\PG(r,\overline{\mathbb{Q}})$, $r \in \mathbb{N} \setminus \{ 0,1\}$.
Then $S$ is not finitely generated.
\eop
\ec

\medskip
\subsection{Infinite dimension}

When $\omega$ is any infinite cardinal number, then for $\PG(\omega,\F)$, with $\F$ either real-closed or algebraically closed, there {\em are} Singer groups.
For,  let $\chi$ be a set of indeterminates such that $\F(\chi)/\F$ has degree $\omega$.
Then as in \S \ref{class}, $\F(\chi)^\times$ acts naturally, linearly and sharply transitively on the vector space $\F^\omega$. Passing to the 
corresponding projective space yields the construction.\\

\bigskip
\vspace*{1cm}\section{Virtual Singer groups and virtual fields}

Let $\Gamma$ be a projective plane, and $Y \leq \Aut(\Gamma)$. We say that $Y$ is {\em virtually Singer} (or $Y$ is a {\em virtual Singer group}) if $Y$ acts 
freely on the points of the plane, and if the number of $Y$-orbits on the point set is finite. (One could also relax this condition by asking that 
the cardinality of $Y$ and of the point set are the same.) We say that $\Gamma$ is {\em virtually Singer} if $\Gamma$ admits a virtual Singer group.
We also say that a field $\K$ is {\em (virtually) Singer} in {\em degree} $m$ ($m \in \mathbb{N}^\times$) if the space $\PG(m,\K)$ is (virtually) Singer. If we do not 
mention the degree, we mean the planar case.\\

\bt[Examples]
We have the following for virtual Singer groups.
\begin{itemize}
\item
All Singer groups are virtually Singer.
\item
All groups acting freely on finite planes are virtually Singer. 
\item
All planes $\PG(2,\F)$ with $\F$ not cubically closed are virtually Singer.
\item
No plane of the statement of Corollary \ref{exmp} can be virtually Singer.
\end{itemize}
\et
{\em Proof}.\quad
All of the statements are straightforward, except the last one. This follows with a similar proof as Corollary \ref{exmp}. 
\eop \\

Two questions which arise are: (1) {\em do there exist real-closed or algebraically closed fields which are virtually Singer ?}; (2) {\em do there exist real-closed or algebraically closed fields which are virtually Singer but not Singer ?}\\

\bigskip
\vspace*{1cm}\section{Singer groups of $\mathbb{F}_1^m$-spaces}

Let $\mathbf{P}$ be an $m$-dimensional projective space over $\mathbb{F}_{1^n}$, where $(n,m) \in \mathbb{N}^\times \times \mathbb{N}^\times$. It is a set 
of $m + 1$ sets $X_i$ of size $n$, each endowed with a free action of the multiplicative group $\mu_n^i \cong C_n$, together with the induced subspace structure.
The linear automorphism group of $\mathbf{P}$  is $C_n \wr \mathbf{S}_{m + 1}$ (elements consist of $(m + 1)\times(m + 1)$-matrices with only one nonzero
entry per row and column, and each such entry is an element of $C_n$). All this can be found in \cite{chap2}.

\subsection{First construction}

It is clear that once we fix a sharply transitive subgroup $S$ of $\mathbf{S}_{m + 1}$, 
we can construct a Singer group $S(n)$ of $\mathbf{P}$ by taking the direct product with a diagonal group $\langle (\nu_1,\nu_2,\ldots,\nu_{m + 1}) \rangle \cong C_n = \langle \nu \rangle$, 
where each $\nu_i$ is a copy of $\nu$ (acting on $X_i$), all $X_j$ being fixed. Now note that 

\begin{equation}
i\ \ \vert\ \ j \ \ \mbox{implies}\ \ S(i) \leq S(j). 
\end{equation}
(For each $S(j)$, we use the same group $S$.) So we obtain a directed system of Singer groups of projective $m$-spaces over all finite extensions of $\mathbb{F}_1$,
and after passing to the limit we obtain a Singer group $S(\infty)$ of $\PG(m,\overline{\Fun})$.

As we fix $S$ in the construction above, passing to the direct limit amounts to constructing $\overline{\Fun}$ by taking unions of all the cyclic groups $C_n$ in the 
appropriate way. As we have seen, there is no way this approach can be pulled to other finite fields (due to the fact that at the $\Fun$-level, we cannot see
the extra needed divisibility condition ``$(j/i,3) = 1$'').

\subsection{General construction}

There is a more generic construction, which in fact captures all possible Singer groups of the spaces $\PG(m,\Fun^d)$.
Let $\mathbf{P}$ be an abstract set (say suggestively of $m + 1$ elements with $m \in \mathbb{N}^\times$), and let $S \leq \mathrm{Sym}(\mathbf{P})$ be a transitive group.
We require that
\begin{itemize}
\item[(CY)] For some element $x \in \mathbf{P}$ (and then all elements), $S_x$ is cyclic. 
\end{itemize}
Let $S_x \cong C_n$; then clearly $S$ has a natural action as a Singer group on $\PG(m,\Fun^n)$, and all Singer groups of this space arise in this way.


\end{document}

%% file: IntroAbsHyper.tex
\section{Absolute Arithmetic}

In a paper which was published in 1957 \cite{anal}, Tits made a seminal and provocative remark which alluded to the fact that through a certain analogy between the groups 
$\mathbf{GL}_n(q)$ (or $\mathbf{PGL}_{n}(q)$), $q$ any prime power, and the symmetric groups $\mathbf{S}_n$, one should interpret $\mathbf{S}_n$ as a Chevalley group ``over the field of characteristic one'':

 \begin{equation}
 \label{GL1}
 \lim_{q \to 1}\mathbf{PGL}_n(q) = \mathbf{S}_n.
 \end{equation}

Only much later serious considerations were made about Tits's point of view, and nowadays a deep theory is being developped on the philosophy over $\mathbb{F}_1$.

In fact, underlying this idea is the fact that thin (spherical) buildings \cite{AbBr,Titslect} are well-defined objects, and with a natural definition of automorphism group, the latter would 
become Weyl groups in thick buildings of the same type defined over ``real fields'' if one considers the appropriate building. So although for instance thin buildings of type $\mathbf{A}_n$ are present in abundance in any 
thick building $\PG(n,q)$ over some finite field $\mathbb{F}_q$, we cannot define them as an incidence geometry over some field, since the cardinality of the latter should be one. Still, the automorphism
group of the underlying geometry (which is just $2^X$ if $X$ is the point set of the thin geometry) would precisely be the Weyl group of the associated thick building | namely the symmetric group on $n$ letters.

Other combinatorial apects of absolute geometry (e.g. concerning projective and affine spaces, and Linear Algebra) were described in an unpublished but important manuscript by Kapranov and Smirnov \cite{KapranovUN} (which we will partly reproduce later on), making also formal sense of the expression (\ref{GL1}).

\medskip
\subsection{Absolute geometry}

The geometry of algebraic curves underlying the structure of global fields of positive characteristic lies at the base of the solution of several deep and fundamental questions in Number Theory. Several formulas of combinatorial nature (such as the number of subspaces of a finite projective space) still keep a meaninful value if evaluated at $q = 1$. Such results seem to suggest the existence of a mathematical object which is a nontrivial limit of finite fields $\mathbb{F}_q$ with $q \rightarrow 1$.
The goal would be to define an analogue, for number fields, of the geometry underlying the arithmetic theory of function fields, cf. \cite{Connes4}.

In the early nineties, Christopher Deninger published his studies (\cite{Deninger1991}, \cite{Deninger1992}, \cite{Deninger1994}) on motives and regularized determinants. In \cite{Deninger1992}, Deninger gave a description of conditions on a category of motives that would admit a translation of Weil's proof of the Riemann Hypothesis for function fields of projective curves over finite fields $\F_q$ to the hypothetical curve $\overline{\Spec(\Z)}$. In particular, he showed that the following formula would hold:

	\begin{eqnarray}
	\label{eq2}
		&2^{-1/2}\pi^{-s/2}\Gamma(\frac s2)\zeta(s) \ &= \\
		&\frac{\det_\infty\Bigl(\frac 1{2\pi}(s-\Theta)\Bigl| H^1(\overline{\Spec(\Z)},\cO_\cT)\Bigr.\Bigr)}{\det_\infty\Bigl(\frac 1{2\pi}(s-\Theta)\Bigl| H^0(\overline{\Spec(\Z)},\cO_\cT)\Bigr.\Bigr)\det_\infty\Bigl(\frac 1{2\pi}(s-\Theta)\Bigl| H^2(\overline{\Spec(\Z)},\cO_\cT)\Bigr.\Bigr)}, &
	\end{eqnarray}
where $\det_\infty$ denotes the regularized determinant, $\Theta$ is an endofunctor that comes with the category of motives and the $H^i(\overline{\Spec(\Z)},\cO_\cT)$ are certain proposed cohomology groups. This description combines with Kurokawa's work on multiple zeta functions (\cite{Kurokawa1992}) from 1992 to  the hope that there are motives $h^0$ (``the absolute point''), $h^1$ and $h^2$ (``the absolute Tate motive'') with zeta functions
	\begin{equation}
	\label{eqzeta}
		\zeta_{h^w}(s) \ = \ \det{}_\infty\Bigl(\frac 1{2\pi}(s-\Theta)\Bigl| H^w(\overline{\Spec\Z},\cO_\cT)\Bigr.\Bigr) 
	\end{equation}
for $w=0,1,2$. Deninger computed that $\zeta_{h^0}(s)=s/2\pi$ and $\zeta_{h^2}(s)=(s-1)/2\pi$. Manin proposed in \cite{Manin} the interpretation of $h^0$ as $\Spec(\Fun)$ and the interpretation of $h^2$ as the affine line over $\Fun$. The search for a proof of the Riemann Hypothesis became a main motivation to look for a geometric theory over $\Fun$. \\

For much more on the Algebraic Geometry in characteristic one, we refer the reader to the monograph \cite{AAbook}.

\medskip
\subsection{Hyperfield extensions and the Krasner hyperfield}

In \cite{Connes3}, Connes and Consani relate hyperstructures to the geometry of $\mathbb{F}_1$, giving rise to an interesting connection between hyperfield extensions of the so-called ``Krasner hyperfield'', and sharply transitive actions (on points) of automorphism groups of  combinatorial projective planes (see e.g. \cite{Hall,Pren,KTDZ} for strongly related discussions). (One is also referred to the paper \cite{Order} by the author for far more details on the Connes-Consani connection, sketched in the framework of  ``prime power conjectures''.) The {\em Krasner hyperfield} $\mathbf{K}$ is the hyperfield $(\{0,1\},+,\cdot)$ with additive neutral element $0$, usual multiplication with identity $ \id$, and satisfying the ``hyperrule''

\begin{equation}
1 + 1 = \{0,1\}.
\end{equation}

In the category of hyperrings, $\mathbf{K}$ can be seen as the natural extension of the commutative pointed monoid $\mathbb{F}_1$, that is, $(\mathbf{K},\cdot) = \mathbb{F}_1$.

 \begin{equation}
  \hspace*{-1cm}
\begin{array}{ccc}
\mathrm{Spec}(\mathbf{K})&\leftarrow&\mathrm{Spec}(\mathbb{Z})\\
&&\\
\downarrow\scriptstyle{}&\swarrow&\\
&&\\
\mathrm{Spec}(\mathbb{F}_1)&&\\


\end{array} 
\end{equation}

\medskip
\subsection{The present paper}

It is precisely a number of foundational questions that arise in \cite{Connes3} in the context of classifying hyperfield extensions of $\mathbf{K}$, and that can be traced back to Hall \cite{Hall} in 1947, which we want to consider in this paper. Formulated in a rather general form, one set of problems we want 
to study is:\\

\quad$\natural$ \textsc{Question(s)}.\quad
{\em For which fields $\K$  can classical planes $\mathbf{PG}(2,\mathbb{K})$ admit sharply transitive automorphism groups?}\\

The precise connection with the field of characteristic one is explained in the next section.
Very roughly, we will see that a group is in some sense ``defined over the field with one element'' if there exists a projective plane such that this group
acts sharply transitively, as an automorphism group, on its points.

And on the other hand, we also want to know\\

\quad$\natural$ \textsc{Question(s)}.\quad
{\em Which groups $G$ act as  sharply transitive automorphism groups of projective planes?}\\

Here, the sharply transitive action is taken on the points. Groups with this action are called {\em Singer groups} (of the planes) throughout.

The reader recalls a basic classical result for finite projective planes which states that if a finite projective plane $\Gamma$ of order $n$ admits the 
projective general linear group $\mathbf{PGL}_3(n)$ as an automorphism group (so $n$ is already assumed to be a prime power), then $\Gamma$
comes from a vector space over $\F_q$, i.e., is coordinatized over $\F_q$, and the group acts doubly transitively on both points and lines. In fact, similar results can be obtained when we ask that the group is abstract and acts doubly transitively on points or lines, or even by merely assuming it is big enough.
In other words, once an automorphism group is assumed to be big or classical enough or {\em acts} classically enough, this can only work for a classical plane, and the group contains  the information of only the classical plane (which can be reconstructed from the group) for this action. So we ask ourselves the 
following\\

\quad$\natural$ \textsc{Question}.\quad
{\em Can one group act as a Singer group on nonisomorphic projective planes?}\\

(Whereas one and the same classical plane could admit nonisomorphic Singer groups.)\footnote{It should be noted that for {\em affine} Singer actions, such as sharply transitive actions on the point set of an affine plane, many examples exist of nonisomorphic planes that admit the same Singer group. Even for buildings of higher rank such affine actions are known.} 

The sharply transitive action on points (or lines) is much more rigid than the ``linear actions'', but in the finite case no examples of nonclassical planes are known which admit sharply transitive automorphism groups. Even when the groups are assumed to be:
\begin{itemize}
\item
abelian, or even
\item
cyclic
\end{itemize}
we do not know that the planes are classical, or even have {\em a prime power order}! 
In fact, we will see that by a result of Karzel, ``abelian'' implies ``cyclic''. And as we will also see, this is not true at all in the infinite case: Karzel showed that infinite finitely generated abelian groups can {\em never} act sharply transitively on classical projective planes, and 
we shall prove that virtually {\em all} infinite abelian groups can act as Singer groups on {\em certain} projective planes, but those will never be classical.

\bt
If an infinite abelian group $A$ does not have involutions, there exists a projective plane $\Gamma$ such that $A \leq \Aut(\Gamma)$, 
and $A$ acts sharply transitively on the points of $\Gamma$.
\et

This result is a corollary of a criterium which naturally generalizes an old result of Hughes from the 1950s. Another corollary of this criterium is:

\bt
For any not necessarily finitely free group $\mathbf{F}$, there exists a projective plane $\Gamma$ such that $\mathbf{F} \leq \Aut(\Gamma)$, 
and $\mathbf{F}$ acts sharply transitively on the points of $\Gamma$.
\et

Hughes obtained the same result for finitely generated free groups (see the next sections for the details).

Eventually, one of the things we want to understand is what the relation is between the structure of a Singer group of a classical plane, and the field over which the plane is  coordinatized. In particular, at the moment we do not know wether, if such a field is algebraically closed or real-closed, such groups exist. For the algebraic closures of {\em finite fields}, however, we will prove indeed that Singer groups cannot exist.

\bt
For any prime $p$, the projective plane $\PG(2,\overline{\F_p})$ does not admit Singer groups. 
\et

It also holds for higher dimensions.
The proof is in some sense related to the Ax-Grothendieck theorem on polynomial maps, and can be situated in the theory of \cite{Serre}. 

As for real-closed fields, for the reals, the same result holds.

\bt
The projective plane $\PG(2,\mathbb{R})$ does not admit Singer groups. 
\et

Again, this is a corollary of a general criterium which leads to the same result for many other real-closed fields (and also for the higher
dimensional problem). In dimension $2$, it says that if a field $\F$ is real-closed, and $\vert \Aut(\F) \vert < \vert \F \vert$, then $\PG(2,\F)$ cannot admit Singer groups.\\

It seems that some general principle is present: the farther away a field $\F$ is from being algebraically closed, the easier it is to construct (isomorphism classes of) Singer groups for projective spaces over $\F$. It would be extremely interesting to find a precise formulation for this principle.

Finally, we want to consider not just planes, but any (also infinite dimensional) projective space for this problem.\\

\subsection{Acknowledgements}

I want to thank Lenny Taelman and Hendrik Lenstra, Jr. for several useful conversations during the Spring 2013 DIAMANT symposium in Heeze. I also
am grateful to Manuel Merida and Lenny Taelman for a number of suggestions on a draft of the present text. Finally, I wish to 
express my gratitude to Adam Topaz and Mariusz Wodzicki for some inspiring conversations during a lecture course I gave on the subject of $\mathbb{F}_1$  at the wonderful Mathematics seminar of UC Berkeley in April 2014.

%% file: F1Hyper.tex
\section{The hyperring of ad\`{e}le classes}

In a recent paper \cite{Connes3}, the authors discovered unexpected connections between hyperrings and (axiomatic) projective geometry, foreseen with certain group actions.
Denoting the profinite completion of $\mathbb{Z}$ by $  \widehat{\mathbb{Z}}$ (and noting that it is isomorphic to the product $\prod_p\mathbb{Z}_p$ of all $p$-adic integer rings), the {\em integral ad\`{e}le ring} is defined as 
\begin{equation}
\mathbb{A}_\mathbb{Z} = \mathbb{R} \times   \widehat{\mathbb{Z}},
 \end{equation}
 
 \noindent
  endowed with the product topology.
Let $\K$ be a global field (that is, a finite extension of $\mathbb{Q}$ or the function field of an algebraic curve over $\mathbb{F}_q$ | the latter is a finite field extension of the field of rational functions $\mathbb{F}_q(X)$). The {\em ad\`{e}le ring} of $\K$ is 
given by the expression
\begin{equation}
\mathbb{A}_{\K} = {\prod_{\nu}^{}}'\K_{\nu},
\end{equation}
which is the restricted product of local fields $\K_{\nu}$, labeled by the places of $\K$. 
If $\K$ is a number field, the ad\`{e}le ring of $\K$ can also be defined as the 
tensor product
\begin{equation}
\mathbb{A}_\K = \K \otimes_{\mathbb{Z}}\mathbb{A}_\mathbb{Z}
\end{equation}
with a suitable topology.

We need a few more definitions.

\subsection{Hyperrings and hyperfields}

Let $H$ be a set, and ``$+$'' be a ``hyperoperation'' on $H$, namely a map

\begin{equation}
+: H \times H \rightarrow (2^H)^\times,
\end{equation}

\noindent
where $(2^H)^\times = 2^H \setminus \{\emptyset\}$. For $U, V \subseteq H$, denote $\{ \cup(u + v)\vert u \in U, v \in V\}$ by $U + V$. (Here, we identify an element $h \in H$ with the singleton $\{h\} \subset H$.)
Then $(H,+)$ is a {\em abelian hypergroup} provided the following properties are satisfied:

\begin{itemize}
\item
$x + y = y + x$ for all $x, y \in H$;
\item
$(x + y) + z = x + (y + z)$ for all $x, y, z \in H$; 
\item
there is an element $0 \in H$ such that $x + 0 = 0 +  x$ for all $x \in H$;
\item
for all $x \in H$ there is a unique $y \in H$ ($=: -x$) such that $0 \in x + y$;
\item
$x \in y + z$ $\Longrightarrow$ $z \in x - y$.
\end{itemize}

\begin{proposition}[\cite{Connes3}]
Let $(G,\cdot)$ be an abelian group, and let $K \subseteq \Aut(G)$. Then the following operation defines a hypergroup structure on $H = \{g^K\vert g \in G\}$:

\begin{equation}
g_1^{K}\cdot g_2^K := (g_1^K\cdot g_2^K)/K.
\end{equation}
\end{proposition}

A {\em hyperring} $(\mathbf{R},+,\cdot)$ is a nonempty set $\mathbf{R}$ endowed with a hyperaddition $+$ and 
a binary operation ``$\cdot$'' which we call ``multiplication'' such that:

\begin{itemize}
\item
$(\mathbf{R},+)$ is an abelian hypergroup with neutral element $0$;
\item
$(\mathbf{R},\cdot)$ is a monoid with multiplicative identity $ \id$;
\item
for all $u, v, w \in \mathbf{R}$ we have that $u\cdot(v + w) = u\cdot v + u\cdot w$ and $(v + w)\cdot u = v\cdot u + w\cdot u$;
\item
for all $u \in \mathbf{R}$ we have that $u\cdot 0 = 0 = 0\cdot u$;
\item
$0 \ne  \id$.
\end{itemize}

A hyperring $(\mathbf{R},+,\cdot)$ is a {\em hyperfield} if $(\mathbf{R} \setminus \{0\},\cdot)$ is a group.

\subsection{The Krasner hyperfield and its extensions}

The {\em Krasner hyperfield} $\mathbf{K}$ is the hyperfield 
\begin{equation}
(\{0,1\},+,\cdot)
\end{equation}
with additive neutral element $0$, usual multiplication with identity $ \id$, and satisfying the hyperrule

\begin{equation}
1 + 1 = \{0,1\}.
\end{equation}

\br[Krasner and $\mathbb{F}_1$]
{\rm
In the category of hyperrings, $\mathbf{K}$ can be seen as the natural extension of the commutative pointed monoid $\mathbb{F}_1$, that is, $(\mathbf{K},\cdot) = \mathbb{F}_1$. 
As remarked in \cite{Connes5}, the Krasner hyperfield {\em encodes the artihmetic of zero and nonzero numbers}, just as $\mathbb{F}_2$ does for even and odd numbers.
(From this viewpoint, it is of no surprise that projective geometry will come into play.)}
\er


Let $\mathbf{R}$ be a commutative ring, and let $G$ be a subgroup of its multiplicative group. The following operations define a hyperring on the set $\mathbf{R}/G$ of $G$-orbits
in $\mathbf{R}$ under multiplication.

\begin{itemize}
\item
\textsc{Hyperaddition}.\quad
$x + y := (xG + yG)/G = \{ xg + yh \vert g,h \in G \}$ for $x, y \in \mathbf{R}/G$.
\item
\textsc{Multiplication}.\quad
$xG\cdot yG := xyG$ for $x, y \in \mathbf{R}/G$.
\end{itemize}

\medskip
\quad\textsc{Important Example}.\quad
Let $\mathbf{R}$ be the finite field $\F_{q^m}$, where $q$ is a prime power and $m \in \mathbb{N}^\times$, and let $G$ be the multiplicative group
$\F_q^{\times} \leq \F_{q^m}^{\times}$. Then we can see $\mathbf{R}$ naturally as an $m$-dimensional $\F_q$-vector space, or better: as an $(m - 1)$-dimensional $\F_q$-projective space. In the latter case, projective points are the cosets $xG$ with $x \ne 0$. And lines, for instance, are of the form $(xG + yG)/G$.
Once one lets $m$ go to $1$, one naturally constructs the Krasner hyperfield $\mathbf{K}$.
These examples will be very important in what is to come. \\

\begin{proposition}[\cite{Connes3}]
Let $\mathbb{K}$ be a field with at least three elements. Then the hyperring $\mathbb{K}/\mathbb{K}^\times$ is isomorphic to the Krasner hyperfield.
If, in general, $\mathbf{R}$ is a commutative ring and $G \subset \mathbb{K}^\times$ is a proper subgroup of the group of units of $\mathbf{R}$, then the hyperring $\mathbf{R}/G$ defined as above contains $\mathbf{K}$ as a subhyperfield if and only if $\{0\} \cup G$ is a subfield of $\mathbf{R}$. 
\end{proposition}

\medskip
\quad\textsc{Important Example [Ad\`{e}le class space and Krasner]}.\quad
Consider a global field $\mathbb{K}$. Its ad\`{e}le class space $\mathbb{H}_\mathbb{K} = \mathbb{A}_{\mathbb{K}}/\mathbb{K}^\times$ is the quotient of a commutative ring $ \mathbb{A}_{\mathbb{K}}$ by $G = \mathbb{K}^\times$, and $\{0\} \cup G = \mathbb{K}$, so it is a hyperring extension of $\mathbf{K}$.\\

\br{\rm
Remark that the ad\`{e}le class space plays a very important role in the non-commutative program of solving the Riemann Hypothesis. (See for instance \cite{Connes1}.)
}
\er

A {\em $\mathbf{K}$-vectorspace} is a hypergroup $E$ provided with a compatible action of $\mathbf{K}$. As $0 \in \mathbf{K}$ acts by retraction (to 
$\{0\} \subset E$) and $ \mathrm{id} \in \mathbf{K}$ acts as the identity on $E$, the $\mathbf{K}$-vectorspace structure is completely determined by the hypergroup structure.
It follows that a hypergroup $E$ is a $\mathbf{K}$-vectorspace if and only if 

\begin{equation}
x + x = \{0,x\}\ \ \mathrm{for}\ \ x \ne 0.
\end{equation}

\medskip
Let $E$ be a $\mathbf{K}$-vectorspace, and define $\mP := E \setminus \{0\}$. For $x, y \ne x \in \mP$, define the {\em line} $L(x,y)$ as 

\begin{equation}
x + y \cup \{x,y\}.  
\end{equation}

\noindent
It can be easily shown | see \cite{Pren} | that $(\mP,\{L(x,y) \vert x, y\ne x \in \mathbf{P}\})$ is a projective space. Conversely, if $(\mP,\mL)$ is the point-line geometry of a projective space with at least $4$ points per line, then a hyperaddition on $E := \mP \cup \{0\}$ can be defined as follows:

\begin{equation}
x + y = xy \setminus \{x,y\}\ \ \mathrm{for}\ \ x \ne y,\ \ \mathrm{and}\ \ x + x =\{0,x\}.
\end{equation}

\medskip
Now let $\mathbb{H}$ be a hyperfield extension of $\mathbf{K}$, and let $(\mP,\mL)$ be the point-line geometry of the associated projective space; then Connes and Consani \cite{Connes3} show that $\mathbb{H}^\times$ induces a so-called ``two-sided incidence group'' (and conversely, starting from such a group $G$, there is a unique hyperfield extension $\mathbb{H}$ of $\mathbf{K}$ such that $\mathbb{H} = G \cup \{0\}$). By the Veblen-Young result, this connection is reflected by the next theorem for the finite case.

\begin{proposition}[\cite{Connes3}]
\label{CCH}
Let $\mathbb{H} \supset \mathbf{K}$ be a finite commutative hyperfield extension of $\mathbf{K}$. Then one of the following cases occurs:
\begin{itemize}
\item[{\rm (i)}]
$\mathbb{H} = \mathbf{K}[G]$ for a finite abelian group $G$.
\item[{\rm (ii)}]
There exists a finite field extension $\mathbb{F}_q \subseteq \mathbb{F}_{q^m}$ such that $\mathbb{H} = \mathbb{F}_{q^m}/\mathbb{F}_q^\times$.
\item[{\rm (iii)}]
There exists a finite nonDesarguesian projective plane admitting a sharply point-transitive automorphism group $G$, and $G$ is the abelian incidence group associated to $\mathbb{H}$. 
\end{itemize}
\end{proposition}


In case (i), there is only one line (otherwise we have to be in the other cases), so for all $x, y, x',  y' \in \mathbb{H} \setminus \{0\}$ with $x \ne y$ and $x' \ne y'$, we must have

\begin{equation}
L(x,y) =  (x + y) \cup \{x,y\} = (x' + y') \cup \{x',y'\} = L(x' + y') = \mathbb{H} \setminus \{0\}. 
\end{equation}

\medskip
\noindent
In other words, hyperaddition is completely determined:

\begin{equation}
\left\{\begin{array}{cc}
x + 0 = x  &\mathrm{for}\ \ x \in \mathbb{H}\\
x + x = \{0,x\}    &\mathrm{for}\ \ x \in \mathbb{H}^\times\\
x + y = \mathbb{H}  \setminus \{0,x,y\}   &\mathrm{for}\ \ x \ne y \in \mathbb{H}^\times\\
\end{array}\right.
\end{equation}

\subsection{The present paper}

There exist infinite hyperfield extensions $\mathbb{H} \supset \mathbf{K}$ for which $\mathbb{H}^\times \cong \mathbb{Z}$ and not coming from Desarguesian projective spaces in the above sense, see M. Hall \cite{Hall}, and the next section. This remark, together with the following general version of 
Theorem \ref{CCH} (see the remark before that theorem), is the starting point of the present paper.

\begin{proposition}[\cite{Connes3}]
Let $\mathbb{H} \supset \mathbf{K}$ be a hyperfield extension of the Krasner hyperfield $\mathbf{K}$. Then there exists a  projective space admitting a sharply point-transitive automorphism group $A$, and $A$ is the incidence group associated to $\mathbb{H}$. 
\end{proposition}

The space could be nonDesarguesian if its dimension is two.
If the dimension of the space is at least three, we know the space {\em is} coordinatized over a skew field by the Veblen-Young result, but when one does not assume the group to be commutative for these spaces, not much seems to be known about such actions. And in the planar case, we can have axiomatic projective planes which are not associated to vector $3$-spaces (over some skew field), and by Hall's result, such planes could admit extremely strange sharply transitive automorphism groups, such as the infinite cyclic group $\mathbb{Z}, +$.

%% file: HyperSinger.bbl
\newpage
{\footnotesize
\begin{thebibliography}{100}
\bibitem{AbBr}
\textsc{P.~Abramenko and K.~S.~Brown}. {\em Buildings. Theory and applications}, Graduate Texts in Math. {\bf 248}, Springer, New York, 2008. 


\bibitem{ArtSchr}
\textsc{E.~Artin and O.~Schreier}.
Eine Kennzeichnung der reell abgeschlossenen K\"{o}rper, {\em Abh. Math. Seminar Univ. Hamburg} {\bf 5} (1927) 225--231.

\bibitem{Baer}
\textsc{R.~Baer}. 
Die Automorphismengruppe eines algebraisch abgeschlossenen K\"{o}rpers der Charakterkistik $0$, {\em Math. Z.} {\bf 117} (1970), 7--17.



\bibitem{Cheretal}
\textsc{G.~Cherlin, T.~Grundh\"{o}fer, A.~Nesin, and H.~V\"{o}lklein}.
Sharply transitive linear groups over algebraically closed fields, 
{\em Proc. Amer. Math. Soc.} {\bf 111} (1991),  541-Ð550.





\bibitem{Connes3}
\textsc{A.~Connes and C.~Consani}.
The hyperring of ad\`{e}le classes, {\em J. Number Theory} {\bf 131} (2011), 159--194.



\bibitem{Connes4}
\bysame.
From monoids to hyperstructures: in search of an absolute arithmetic, in {\em Casimir Force, Casimir Operators and the Riemann Hypothesis}, de Gruyter (2010), pp. 147--198.


\bibitem{Connes2009}
\bysame. 
Schemes over $\F_1$ and zeta functions, {\em Compositio Math.} {\bf 146} (2010), 1383--1415.


\bibitem{Connes5}
\bysame.
On the notion of geometry over $\mathbb{F}_1$, {\em J. Algebraic Geom.}, To appear.





\bibitem{Connes1}
\textsc{A.~Connes, C.~ Consani and M.~ Marcolli}.
 The Weil proof and the geometry of the ad\`{e}les class space, in {\em Algebra, Arithmetic, and Geometry: in honor of Yu. I. Manin, Vol. I}, Progr. Math. {\bf 269}, 2009, Birkh\"{a}user Boston, Inc., Boston, MA, 
pp. 339--405.

\bibitem{Connes2}
 \bysame.
  Fun with $\Bbb F_1$, {\em J. Number Theory} {\bf 129} (2009),  1532--1561.


\bibitem{HSMC}
\textsc{H.~S.~M.~Coxeter}. 
The complete enumeration of finite groups of the form $R_i^2 = (R_iR_j)^{k_{ij}} = 1$, {\em J. London Math. Soc.} {\bf 10} (1935), 21--25.




\bibitem{Deitmarschemes2}
\textsc{A.~Deitmar}. 
Schemes over $\mathbb{F}_1$, in {\em Number Fields and Function Fields Ñ Two Parallel Worlds}, pp. 87Ð100, Progr. Math. {\bf 239}, Birkh\"{a}user Boston, Boston, MA, 2005. 

\bibitem{Deitmarschemes1}
\bysame.
Remarks on zeta functions and $K$-theory over $\mathbb{F}_1$, {\em Proc. Japan Acad. Ser. A Math. Sci.} {\bf 82} (2006),  141--146.

\bibitem{Deitmartoric}
\textsc{A.~Deitmar}.
$\mathbb{F}_1$-Schemes and toric varieties,
{\em Beitr\"{a}ge Algebra Geom.} {\bf 49} (2008), 517--525. 


\bibitem{Deninger1991}
\textsc{C.~Deninger.}
\newblock On the {$\Gamma$}-factors attached to motives,
\newblock {\em Invent. Math.}  {\bf 104} (1991),  245--261.

\bibitem{Deninger1992}
\bysame.
\newblock Local {$L$}-factors of motives and regularized determinants,
\newblock {\em Invent. Math.} {\bf 107} (1992), 135--150.

\bibitem{Deninger1994}
\bysame.
\newblock Motivic {$L$}-functions and regularized determinants,
\newblock in {\em Motives ({S}eattle, {WA}, 1991)}, Proc. Sympos. Pure Math. {\bf 55}, 1994, pp. 707--743.




\bibitem{Hall}
\textsc{M.~Hall, Jr}. Cyclic projective planes, {\em Duke Math. J.} {\bf 14} (1947), 1079--1090.

\bibitem{Hughes}
\textsc{D.~R.~Hughes}. A note on difference sets, {\em Proc. Amer. Math. Soc.} {\bf 6} (1955), 689--692. 

\bibitem{KapranovUN}
\textsc{M.~Kapranov and A.~Smirnov.}
\newblock Cohomology determinants and reciprocity laws: number field case,
\newblock  Unpublished preprint.

\bibitem{Karzel}
\textsc{H.~Karzel}.
Projektive R\"{a}ume mit einer kommutativen transitiven Kollineationsgruppe,
{\em Math. Z.} {\bf 87} 1965 74--77. 


\bibitem{8}
\textsc{B.~Kleiner and B.~ Leeb}. 
Rigidity of quasi-isometries for symmetric spaces and Euclidean buildings, {\em Inst. Hautes \'{E}tudes Sci. Publ. Math.} {\bf 86} (1997), 115--197.

\bibitem{Kurokawa1992}
\textsc{N.~Kurokawa.}
\newblock Multiple zeta functions: an example,
\newblock in {\em Zeta Functions in Geometry ({T}okyo, 1990)}, Adv. Stud. Pure Math. {\bf 21}, pp. 219--226, 1992.

\bibitem{Kurokawa2005}
\bysame.
\newblock Zeta functions over $\Fun$,
\newblock {\em Proc. Japan Acad. Ser. A Math. Sci.} {\bf 81} (2005), 180--184. 


\bibitem{Kurokawa2002}
\textsc{N.~Kurokawa, H.~Ochiai, M.~Wakayama.}
\newblock Absolute derivations and zeta functions,
\newblock {\em Doc. Math.}, Extra Volume: Kazuya Kato's Fiftieth Birthday (2003), 565--584.


\bibitem{LL}
\textsc{J.~L\'opez Pe\~na and O.~Lorscheid.}
\newblock Torified varieties and their geometries over $\Fun$,
\newblock  {\em Math. Z.} {\bf 267} (2011), 605--643.

\bibitem{LLb}
\bysame.
Mapping $\Fun$-land: An overview of geometries over the field with one element, 
 in {\em Noncommutative Geometry, Arithmetic, and Related Topics}, pp. 241--265, Johns Hopkins Univ. Press, Baltimore, MD, 2011.

\bibitem{Lorscheid2009}
\textsc{O.~Lorscheid.}
\newblock Algebraic groups over the field with one element,
{\em Math. Z.} {\bf 271} (2012), 117--138.



\bibitem{Manin}
\textsc{Yu.~Manin}.
 Lectures on zeta functions and motives (according to Deninger and Kurokawa), Columbia University Number Theory Seminar (New York, 1992), {\em Ast\'{e}risque} {\bf 228} (1995), 121--163.
 

\bibitem{Manin2008}
\textsc{Y.~Manin}.
\newblock Cyclotomy and analytic geometry over $\Fun$, in
{\em Quanta of Maths}, pp. 385--408, Clay Math. Proc. {\bf 11}, Amer. Math. Soc., Providence, RI, 2010.



\bibitem{Pren}
\textsc{W.~Prenowitz}. 
Projective geometries as multigroups, {\em Amer. J. Math.} {\bf 65} (1943), 235--256.




\bibitem{Serre}
 \textsc{J.-P.~Serre}. How to use finite fields for problems concerning infinite fields, in {\em Arithmetic, Geometry, Cryptography and Coding Theory}, Contemp. Math. {\bf 487}, Providence, R.I.: Amer. Math. Soc., pp. 183Ð193.

\bibitem{Soule}
\textsc{C.~Soul\'{e}}. 
Les variŽt\'{e}s sur le corps \`{a} un \'{e}l\'{e}ment, {\em Mosc. Math. J.} {\bf 4} (2004), 217--244, 312.


\bibitem{Soule2009}
\textsc{C.~Soul\'{e}}.
\newblock Lectures on algebraic varieties over $\Fun$, Preprint.


\bibitem{Order}
\textsc{K.~Thas}. Order in building theory, in {\em Surveys in Combinatorics 2011}, London Math. Soc. Lecture Notes Series {\bf 392}, Cambridge Univ. Press, pp. 235--331, 2011.

\bibitem{chap2}
\textsc{K.~Thas}.
{\em The Weyl-functor | Introduction to absolute arithmetic}, in  {\em Absolute Arithmetic and $\mathbb{F}_1$-Geometry} (ed. K. Thas), 
With chapters by Borger, Deitmar, Le Bruyn, Lorscheid, Manin, Thas, EMS Publishing House, 38pp., submitted.


\bibitem{chap}
\textsc{K.~Thas}.
{\em The combinatorial-motivic nature of $\mathbb{F}_1$-schemes}, in  {\em Absolute Arithmetic and $\mathbb{F}_1$-Geometry} (ed. K. Thas), 
With chapters by Borger, Deitmar, Le Bruyn, Lorscheid, Manin, Thas, EMS Publishing House, 76pp., submitted.

\bibitem{AAbook}
\textsc{K.~Thas}. (Ed.)
{\em Absolute Arithmetic and $\mathbb{F}_1$-Geometry}, 
With chapters by Borger, Deitmar, Le Bruyn, Lorscheid, Manin, Thas, EMS Publishing House, 378pp., submitted.



 

\bibitem{KTDZ}
\textsc{K.~Thas and D.~ Zagier}.
Finite projective planes, Fermat curves, and Gaussian periods,
{\em J. European Math. Soc. (JEMS)} {\bf 10} (2008), 173--190. 

\bibitem{anal}
\textsc{J.~Tits}.
Sur les analogues alg\'{e}briques des groupes semi-simples complexes,
{\em Centre Belge Rech. Math., Colloque d'Alg\`{e}bre Sup\'{e}rieure, Bruxelles du 19 au 22 d\'{e}c. 1956}, pp. 261-289, 1957.




\bibitem{Tits}
\textsc{J.~Tits}. 
Sur la trialit\'{e} et certains groupes qui s'en d\'{e}duisent,
{\em Inst. Hautes Etudes Sci. Publ. Math.} {\bf 2}  (1959), 13--60. 


\bibitem{Titslect}
\bysame.
{\em Buildings of Spherical Type and Finite BN-Pairs},
Lecture Notes in Mathematics {\bf 386},
Springer-Verlag, Berlin|New York, 1974.
 













\bibitem{TitsWeiss}
\textsc{J.~Tits and R.~M.~ Weiss}.
{\em Moufang Polygons}, Springer Monographs in Mathematics, Springer-Verlag, Berlin, 2002.

\bibitem{Toen2008}
\textsc{B.~To\"en and M.~Vaqui\'e}.
\newblock Au-dessous de $Spec\; \mathbb{Z}$,
\newblock {\em J. K-Theory} {\bf 3} (2009), 1--64.

\bibitem{POL}
\textsc{H.~ Van Maldeghem}.
{\em Generalized Polygons}, Monographs in Mathematics {\bf 93}, Birkh\"{a}user-Verlag, Basel, 1998.



\end{thebibliography}}
